\documentclass[11pt]{amsart}
\usepackage{geometry}
\usepackage{amsmath,amsthm,amssymb}
\usepackage{graphicx}
\usepackage{latexsym}
\usepackage{tikz}
\usetikzlibrary{arrows.meta,calc}
\usepackage[hidelinks]{hyperref}

\usepackage{amsfonts}

\numberwithin{equation}{section}

\theoremstyle{plain}
\newtheorem{theorem}{Theorem}[section]
\newtheorem{lemma}[theorem]{Lemma}
\newtheorem{corollary}[theorem]{Corollary}
\newtheorem{proposition}[theorem]{Proposition}

\theoremstyle{definition}

\newtheorem{remark}[theorem]{Remark}

\begin{document}

\title[Subharmonicity, convexity and Fourier transform]{Subharmonic functions and convexity of the Radon\\ transform,
with applications to the Fourier transform}

\author{A. Iosevich}
\thanks{The first author was supported in part by NSF grant DMS-2506858 and a 2026 Simons Fellowship.}
\address{Department of Mathematics, University of Rochester, NY, USA}
\email{iosevich@gmail.com}

\author{E. Liflyand}
\address{Department of Mathematics, Bar-Ilan University, 52900 Ramat-Gan, Israel}
\email{liflyand@gmail.com}



\date{}

\begin{abstract}
Integrating a subharmonic function over parallel hyperplanes produces a convex function of the signed distance from the origin, provided that the function is absolutely integrable on bounded slabs. We prove this first for twice continuously differentiable functions and then without any smoothness assumption, and the same argument applies to a function defined only on a slab. Examples on the disk explain why subharmonicity on a bounded convex domain does not imply convexity of the integrals over its chords.

We then apply a one-dimensional Fourier representation for convex functions to the Radon transform, which leads to a directional formula for integrable subharmonic functions on slabs. For functions on bounded convex domains whose Radon transforms are convex, we obtain a weighted estimate for the integral of the Fourier transform in terms of boundary layers and the gradient.
\end{abstract}

\maketitle
\enlargethispage{2pt}

\section{Introduction}

The multidimensional Fourier transform
\begin{equation}\label{gft}
\widehat{f}(x)=\int_{\mathbb R^d}f(u)e^{-i\langle x,u\rangle}du,
\end{equation}
where $\langle x,u\rangle=x_1u_1+\cdots+x_du_d$, reduces to the one-dimensional Fourier transform of the Radon transform
\begin{equation}\label{t5.1}
Rf(\theta,t)=\int_{\langle x,\theta\rangle=t}f(x)d\omega(x).
\end{equation}
Here $\theta$ belongs to the unit sphere $\mathbb S^{d-1}$, $d\omega$ is the $(d-1)$-dimensional Lebesgue measure on the hyperplane, and $t$ is the signed distance of the hyperplane from the origin. For $f\in L^1(\mathbb R^d)$ and each fixed $\theta$, the Radon transform is defined for almost every $t$ and belongs to $L^1(\mathbb R)$.

The two transforms are related by the Fourier slice theorem, a consequence of Fubini's theorem, which states that
\begin{equation}\label{t5.9}
\widehat f(t\theta)=\int_{-\infty}^{\infty}e^{-its}Rf(\theta,s)ds.
\end{equation}
We refer to \cite{Helgason,Natterer} for the Radon transform and to \cite{IL} for its use in Fourier estimates. The slice theorem is a tempting device, since it turns any one-dimensional result about the Fourier transform into a multidimensional one, but the price is that the one-dimensional hypotheses land on $Rf$, whereas one would much rather impose them on $f$ itself. The purpose of this paper is to describe a situation in which a property of $f$, namely subharmonicity, is transmitted to $Rf$ as convexity, and convexity is exactly the kind of property that one-dimensional Fourier estimates are built on.

The relation between subharmonicity and convexity of the Radon transform was considered in \cite{L2025}, and the calculation behind it is simple, since the second derivative of a hyperplane integral is the integral of the Laplacian, provided that a tangential boundary term vanishes. What requires care is the domain on which the function is subharmonic, because integration over an entire hyperplane and integration over a chord of a bounded domain behave quite differently. Even a constant function on a disk has strictly concave chord integrals.

Our first result makes the argument over whole hyperplanes precise under a domination assumption on the first two derivatives. We then show that the smoothness is unnecessary, since absolute integrability on bounded slabs already produces a unique finite convex representative of $Rf(\theta,\cdot)$. The proof is local in the variable $t$, so it applies equally well when $f$ is defined only on a slab. These results, together with the examples that explain their assumptions, occupy the next three sections.

In Section~\ref{odres} we turn to the one-dimensional Fourier transform, give the representation on which everything else rests, and prove the integrable remainder estimate for convex functions. Section~\ref{gerai} begins with the corresponding directional formula, which we apply first to subharmonic functions on slabs and then to a weighted Fourier estimate on bounded convex domains. In the latter result the convexity of the Radon transform is an explicit assumption, since subharmonicity inside the domain does not imply it. In the Sobolev class used for this application, a stronger weighted integral estimate is classical, and we give the comparison in Section~\ref{gerai}. The last subsection explains what this means for convex functions on bounded domains and why there is no nonzero convex example on a slab.

We write $C_c^\infty$ for the class of smooth functions with compact support. Throughout the paper $C$ denotes a positive constant whose value may change from one occurrence to the next, and $C_D$ indicates that the constant depends on a domain $D$.

\section{Where the convexity comes from}\label{sec_smooth}

Let us begin with a smooth function $u$ in the plane and consider $g(t)=\int_{\mathbb R}u(t,y)dy$. We would like to show that $g$ is convex, and the natural thing to do is to differentiate twice. If we may differentiate under the integral sign, then
$$
g''(t)=\int_{\mathbb{R}}u_{tt}(t,y)dy=\int_{\mathbb{R}}\Delta u(t,y)dy-\int_{\mathbb{R}}u_{yy}(t,y)dy,
$$
where $\Delta u=u_{tt}+u_{yy}$ is the Laplacian. If $u$ is subharmonic, which for a smooth function means that $\Delta u\geq0$, the first integral on the right is nonnegative, while the second is a boundary term, and if it vanishes we are done. The whole question therefore comes down to making sure that the integrals exist and that nothing is left over at infinity.

The same calculation works in every dimension. Let $d\geq2$, fix $\theta\in\mathbb S^{d-1}$, let $\theta^\perp=\{y\in\mathbb R^d\mid\langle y,\theta\rangle=0\}$, and write $x=t\theta+y$ with $t\in\mathbb{R}$ and $y\in\theta^\perp$. In these coordinates the Radon transform in the direction $\theta$ takes the form
\begin{equation}\label{eq_radon}
Ru(\theta,t)=\int_{\theta^\perp}u(t\theta+y)dy,
\end{equation}
where $dy$ is the $(d-1)$-dimensional Lebesgue measure on $\theta^\perp$. We write $\Delta_y$ for the Laplacian in the variable $y\in\theta^\perp$, so that $\Delta u=\partial_t^2u+\Delta_yu$. The direction $\theta$ remains fixed throughout, and convexity always refers to the variable $t$.

\begin{proposition}\label{prop_smooth}
Let $u\in C^2(\mathbb{R}^d)$ satisfy $\Delta u\geq0$, and fix $\theta\in\mathbb S^{d-1}$. Suppose that for every compact interval $J\subset\mathbb{R}$ there is a nonnegative function $H_J\in L^1(\theta^\perp)$ such that
\begin{equation}\label{eq_domination}
|D^\alpha u(t\theta+y)|\leq H_J(y)
\end{equation}
for all $t\in J$, $y\in\theta^\perp$ and $|\alpha|\leq2$, where $D^\alpha$ denotes a coordinate derivative in the orthonormal coordinates $(t,y)$ and $D^0u=u$. Then $Ru(\theta,\cdot)$ is a finite convex function of class $C^2$, and more precisely
\begin{equation}\label{eq_identity}
\frac{d^2}{dt^2}Ru(\theta,t)=\int_{\theta^\perp}\Delta u(t\theta+y)dy\geq0.
\end{equation}
\end{proposition}

\begin{proof}
The domination assumption allows us to differentiate twice under the integral sign, and dominated convergence shows that the resulting second derivative is continuous. Since $\Delta u=\partial_t^2u+\Delta_yu$,
\begin{equation}\label{eq_decomposition}
\frac{d^2}{dt^2}Ru(\theta,t)=\int_{\theta^\perp}\Delta u(t\theta+y)dy-\int_{\theta^\perp}\Delta_yu(t\theta+y)dy,
\end{equation}
where both integrals converge absolutely by \eqref{eq_domination}, and it remains to check that the second one vanishes. A cutoff takes care of this. Choose $\chi\in C_c^\infty(\theta^\perp)$ with $0\leq\chi\leq1$, equal to one on the unit ball and to zero outside the ball of radius $2$, and for $R>0$ put $\chi_R(y)=\chi(\frac{y}{R})$. Then $\chi_R(y)\to1$ for every $y$ as $R\to\infty$, and $|\Delta_y\chi_R(y)|\leq\frac{C}{R^2}$ with $C$ independent of $R$. Since $\chi_R$ has compact support, two integrations by parts produce no boundary contribution, and
$$
\int_{\theta^\perp}\chi_R(y)\Delta_yu(t\theta+y)dy=\int_{\theta^\perp}u(t\theta+y)\Delta_y\chi_R(y)dy.
$$
The right side is at most $\frac{C}{R^2}\int_{\theta^\perp}|u(t\theta+y)|dy$ in absolute value, which tends to zero as $R\to\infty$, while the left side tends to $\int_{\theta^\perp}\Delta_yu(t\theta+y)dy$ by dominated convergence and \eqref{eq_domination}. This integral is therefore zero, the equation \eqref{eq_decomposition} reduces to \eqref{eq_identity}, and since the second derivative is nonnegative, the Radon transform is a convex function of $t$.
\end{proof}

\begin{figure}[htbp]
\centering
\begin{tikzpicture}[x=1cm,y=1cm,font=\small,>={Stealth}]
 \fill[blue!7] (-4,-2.4) rectangle (4,2.4);
 \fill[blue!17] (-4,-1.2) rectangle (4,1.2);
 \fill[orange!20,opacity=0.7] (-1,-2.9) rectangle (1.6,2.9);
 \draw[->] (-4.3,0)--(4.6,0) node[right] {$t$};
 \draw[->] (0,-3)--(0,3.1) node[above] {$\theta^\perp$};
 \foreach \t in {-0.6,0.5,1.2} {\draw[thick] (\t,-2.9)--(\t,2.9);}
 \node[above] at (1.4,2.9) {$t\theta+\theta^\perp$};
 \node[below] at (-1,0) {$a$};
 \node[below] at (1.6,0) {$b$};
 \draw (4,1.2)--(4.15,1.2) node[right] {$R$};
 \draw (4,2.4)--(4.15,2.4) node[right] {$2R$};
 \node[anchor=east] at (-4.1,0.6) {$\chi_R=1$};
 \node[anchor=east] at (-4.1,1.8) {$|\Delta_y\chi_R|\leq\frac{C}{R^2}$};
 \draw[->,very thick] (-3.6,-2.7)--(-2.7,-2.7) node[right] {$\theta$};
\end{tikzpicture}
\caption{The geometry of the proofs in the plane, with $\theta=(1,0)$. The Radon transform integrates $u$ over the lines $t\theta+\theta^\perp$. The cutoff $\chi_R$ equals one in the darker band, its second derivatives, of size at most $\frac{C}{R^2}$, are supported in the lighter one, and the proof of Theorem~\ref{thm_weak} only uses the integrability of $u$ over slabs $a<t<b$.}
\label{fig:slab}
\end{figure}
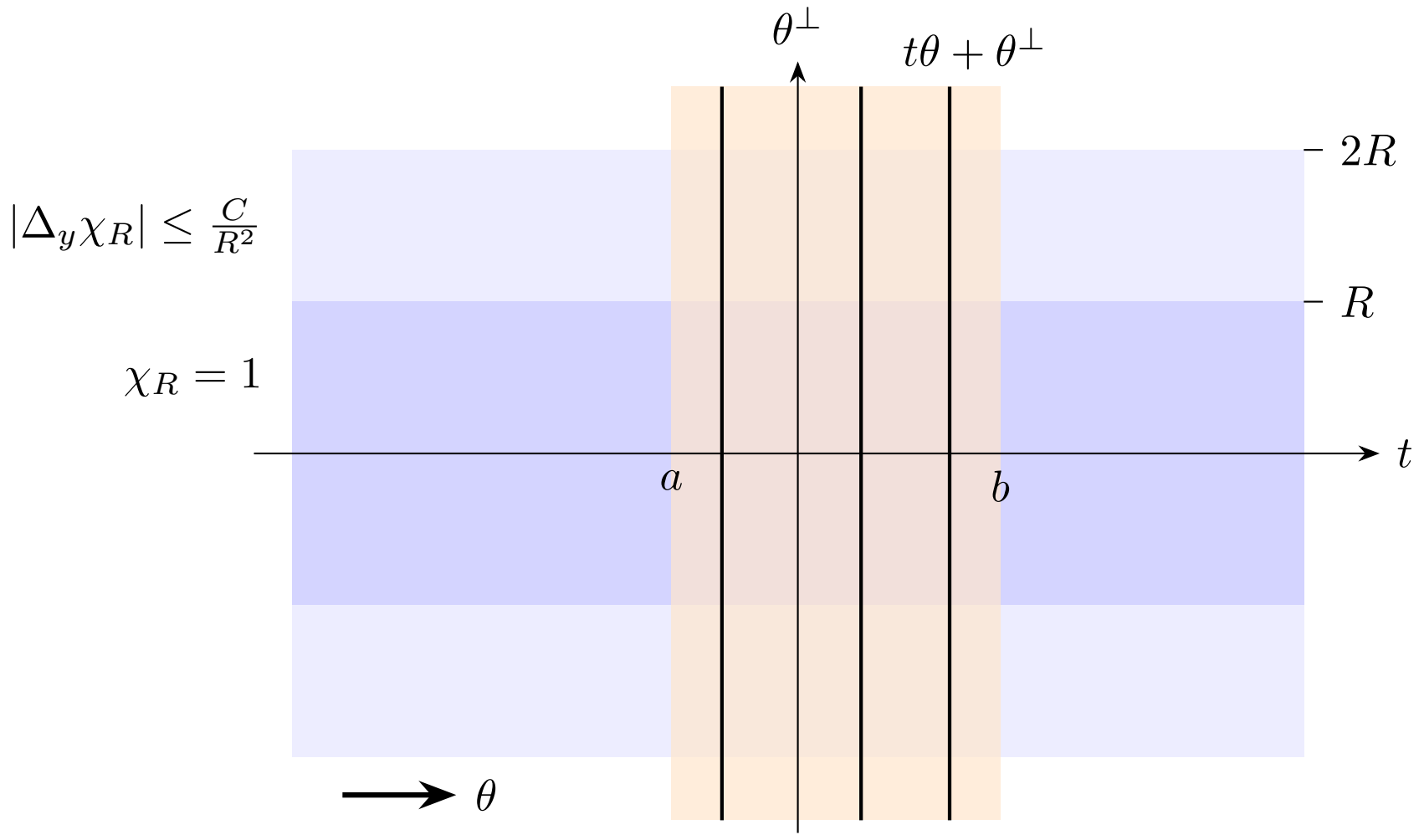

The identity \eqref{eq_identity} is the classical relation between the Radon transform and the Laplacian, according to which the Radon transform intertwines $\Delta$ with $\frac{d^2}{dt^2}$. It is usually proved for rapidly decreasing functions, as in \cite{Helgason,Natterer}, whereas the argument above only needs the domination hypothesis \eqref{eq_domination}, which allows functions that grow exponentially in the direction $\theta$, such as the first example in Section~\ref{sec_examples}. Figure~\ref{fig:slab} shows the geometry of the argument.

\section{Removing the smoothness assumption}\label{sec_weak}

The proof we have just given used more regularity than the conclusion requires. Both derivatives can be moved onto a test function, and once this is done the same argument goes through with no smoothness at all.

Recall that a subharmonic function $u$ on $\mathbb{R}^d$, not identically $-\infty$, is locally integrable, that is, it belongs to $L^1_{\mathrm{loc}}(\mathbb{R}^d)$, and that it satisfies
\begin{equation}\label{eq_subharmonic}
\int_{\mathbb{R}^d}u(x)\Delta\psi(x)dx\geq0
\end{equation}
for every nonnegative $\psi\in C_c^\infty(\mathbb{R}^d)$. This is what it means for $\Delta u\geq0$ to hold in the sense of distributions \cite{ArmitageGardiner,Hormander}, and \eqref{eq_subharmonic} is the only property of $u$ that we shall use.

\begin{theorem}\label{thm_weak}
Let $u$ be subharmonic on $\mathbb{R}^d$, fix $\theta\in\mathbb S^{d-1}$, and suppose that
\begin{equation}\label{eq_slab}
\int_a^b\int_{\theta^\perp}|u(t\theta+y)|dy dt<\infty
\end{equation}
for every bounded interval $(a,b)$. Then $Ru(\theta,t)$ is defined and finite for almost every $t$, and it agrees almost everywhere with a unique finite convex function on $\mathbb{R}$.

If the integral in \eqref{eq_radon} defines a continuous function of $t$ on $\mathbb{R}$, that function is convex everywhere.
\end{theorem}

\begin{proof}
By Fubini's theorem and \eqref{eq_slab}, the function
$$
g(t)=\int_{\theta^\perp}u(t\theta+y)dy
$$
is defined almost everywhere and belongs to $L^1_{\mathrm{loc}}(\mathbb{R})$. We shall show that
\begin{equation}\label{eq_convex_test}
\int_{\mathbb{R}}g(t)\varphi''(t)dt\geq0
\end{equation}
for every nonnegative $\varphi\in C_c^\infty(\mathbb{R})$, which is the distributional form of convexity, and then explain why a function with this property has a convex representative.

With the cutoffs of the preceding proof we set $\psi_R(t\theta+y)=\varphi(t)\chi_R(y)$. This function is smooth, nonnegative and compactly supported, and since
$$
\Delta\psi_R(t\theta+y)=\varphi''(t)\chi_R(y)+\varphi(t)\Delta_y\chi_R(y),
$$
the inequality \eqref{eq_subharmonic} takes the form
\begin{equation}\label{eq_cutoff}
0\leq I_1(R)+I_2(R),
\end{equation}
where
$$
I_1(R)=\int_{\mathbb{R}}\int_{\theta^\perp}u(t\theta+y)\varphi''(t)\chi_R(y)dy dt
$$
and
$$
I_2(R)=\int_{\mathbb{R}}\int_{\theta^\perp}u(t\theta+y)\varphi(t)\Delta_y\chi_R(y)dy dt.
$$
Choose a bounded interval $(a,b)$ containing the support of $\varphi$. Then
$$
|I_2(R)|\leq\frac{C\|\varphi\|_\infty}{R^2}\int_a^b\int_{\theta^\perp}|u(t\theta+y)|dy dt,
$$
which tends to zero, while dominated convergence and Fubini's theorem show that $I_1(R)$ tends to $\int_{\mathbb R}g(t)\varphi''(t)dt$. Letting $R\to\infty$ in \eqref{eq_cutoff}, we arrive at \eqref{eq_convex_test}.

It remains to pass from \eqref{eq_convex_test} to a convex representative. This is the elementary characterization of convexity by a nonnegative second distributional derivative \cite{Hormander}, and we include the proof for the sake of completeness. Choose a nonnegative smooth function $\eta$ supported in $[-1,1]$ with $\int_{\mathbb R}\eta=1$, set $\eta_\varepsilon(t)=\varepsilon^{-1}\eta(\frac{t}{\varepsilon})$, and put
$$
g_\varepsilon(t)=(g*\eta_\varepsilon)(t)=\int_{\mathbb R}g(s)\eta_\varepsilon(t-s)ds.
$$
The compact support of $\eta_\varepsilon$ makes this integral finite, and differentiation under the integral sign shows that $g_\varepsilon$ is smooth. For nonnegative $\varphi\in C_c^\infty(\mathbb{R})$, the function $\widetilde\eta_\varepsilon*\varphi$, where $\widetilde\eta_\varepsilon(s)=\eta_\varepsilon(-s)$, is again smooth, nonnegative and compactly supported, so by \eqref{eq_convex_test}
$$
\int_{\mathbb{R}}g_\varepsilon''\varphi dt=\int_{\mathbb{R}}g_\varepsilon\varphi''dt=\int_{\mathbb{R}}g(\widetilde\eta_\varepsilon*\varphi)''dt\geq0.
$$
If the continuous function $g_\varepsilon''$ were negative at some point, a nonnegative test function supported in a small enough neighborhood of that point would contradict this inequality. Hence $g_\varepsilon''\geq0$, and each $g_\varepsilon$ is convex.

The functions $g_\varepsilon$ converge to $g$ in $L^1_{\mathrm{loc}}(\mathbb R)$ as $\varepsilon\to0$, and we would like to pass convexity to the limit. Convergence in $L^1$ does not provide pointwise convergence everywhere, and it is convexity itself that supplies the missing control.

Let $h$ be a convex function. Averaging the inequality $2h(x)\leq h(x-v)+h(x+v)$ over $0\leq v\leq1$, we get
$$
h(x)\leq\frac{1}{2}\int_{x-1}^{x+1}h(s)ds\leq\frac{1}{2}\int_{x-1}^{x+1}|h(s)|ds,
$$
so a bound for the integral of $|h|$ over a larger interval bounds $h$ from above on a smaller one. For a lower bound, suppose that $h\leq K$ on an interval $I$ of length at least two, where $K\geq0$, and take $x_0\in I$. By convexity,
$$
h\left(\frac{x_0+z}{2}\right)\leq\frac{h(x_0)+h(z)}{2}\leq\frac{h(x_0)+K}{2}
$$
for every $z\in I$. If $h(x_0)+K\geq0$, then $h(x_0)\geq-K$. Otherwise $h$ is at most the negative number $\frac{h(x_0)+K}{2}$ on the set of midpoints $\frac{x_0+z}{2}$, an interval of length at least one contained in $I$, and integrating $|h|$ over this interval we find that $h(x_0)\geq-K-2\int_I|h(s)|ds$. Given a compact interval $I$ of length at least two, the first estimate bounds $g_\varepsilon$ from above on $I$ by its $L^1$ norm over the interval enlarged by one unit at each end, and the second then bounds $g_\varepsilon$ from below on $I$. Since the local $L^1$ norms of the functions $g_\varepsilon$ are uniformly bounded for $0<\varepsilon\leq1$, we conclude that the functions $g_\varepsilon$ are uniformly bounded on each compact interval.

The slopes are controlled in the same way. If $|h|\leq M$ on $[A-1,B+1]$ and $A\leq x<y\leq B$, then by convexity
$$
\frac{h(x)-h(A-1)}{x-(A-1)}\leq\frac{h(y)-h(x)}{y-x}\leq\frac{h(B+1)-h(y)}{B+1-y},
$$
and the outer quotients are bounded in absolute value by $2M$, since their denominators are at least one. The functions $g_\varepsilon$ are therefore uniformly Lipschitz on each compact interval. We apply the Arzel\`a--Ascoli theorem to the functions $g_{\frac{1}{n}}$ on $[-1,1]$, then to a subsequence on $[-2,2]$, and so on, and a diagonal subsequence converges uniformly on every compact interval to a finite function $G$. Passing to the limit in the convexity inequality, we see that $G$ is convex, and since the same subsequence converges to $g$ in $L^1$ on compact intervals, $G=g$ almost everywhere.

A finite convex function on $\mathbb{R}$ is continuous, so two such functions which agree almost everywhere agree everywhere, and this proves the uniqueness. The same observation proves the last assertion when the original Radon integral is continuous.
\end{proof}

The proof is local in the variable $t$, which allows us to work on a slab instead of the whole space.

\begin{corollary}\label{cor_slab}
Let $\theta\in\mathbb S^{d-1}$ and $-\infty<a<b<\infty$. Suppose that $u$ is subharmonic on $D=\{x\in\mathbb R^d\mid a<\langle x,\theta\rangle<b\}$ and that $\int_J\int_{\theta^\perp}|u(t\theta+y)|dy dt<\infty$ for every compact interval $J\subset(a,b)$. Then $Ru(\theta,\cdot)$ agrees almost everywhere on $(a,b)$ with a unique finite convex function. If the Radon integral is continuous on $(a,b)$, it is convex there.
\end{corollary}

\begin{proof}
We use test functions $\varphi\in C_c^\infty(a,b)$ in the proof of Theorem~\ref{thm_weak}, so that each function $\varphi(t)\chi_R(y)$ is supported in $D$. The mollification is performed only at points $t$ with $[t-\varepsilon,t+\varepsilon]\subset(a,b)$, and the compactness argument is applied on compact subintervals of $(a,b)$, so the proof goes through without any assumption at the two endpoints.
\end{proof}

Theorem~\ref{thm_weak} is the counterpart, for parallel hyperplanes, of the classical convexity theorem for spherical means. According to that theorem, the average of a subharmonic function over the sphere of radius $r$ about a fixed point is a convex function of $\log r$ in the plane and of $r^{2-d}$ when $d\geq3$ \cite{ArmitageGardiner,HaymanKennedy}. In both cases the variable is the one in which the relevant harmonic functions are affine, that is, $t$ for functions depending only on $\langle x,\theta\rangle$, and $\log r$ or $r^{2-d}$ for radial functions.

\section{Checking the assumptions on examples}\label{sec_examples}

The integrability assumption does allow nonconstant subharmonic functions. In the plane, consider $u(x,y)=e^{\sqrt{2}x-y^2}$. A direct calculation shows that $u_{xx}=2u$ and $u_{yy}=(4y^2-2)u$, so that $\Delta u(x,y)=4y^2e^{\sqrt{2}x-y^2}\geq0$, and for $\theta=(1,0)$ we find that $Ru(\theta,t)=\int_{\mathbb{R}}e^{\sqrt{2}t-y^2}dy=\sqrt{\pi}e^{\sqrt{2}t}$, a convex function, as expected. For $t$ in a compact interval $J$, the function $u$ and its derivatives of order at most two are bounded by $C_J(1+y^2)e^{-y^2}$, so the domination hypothesis of Proposition~\ref{prop_smooth} holds as well, even though $u$ grows exponentially in the direction $\theta$. The assumptions concern the chosen direction only, and in this example the integrals over the horizontal lines, which correspond to $\theta=(0,1)$, are infinite.

On the other hand, an arbitrary subharmonic function need not have a finite Radon transform at all. The function $u(x)=|x|^2$ is subharmonic, but its integral over every affine hyperplane is infinite.

\begin{figure}[htbp]
\centering
\begin{tikzpicture}[font=\small,>={Stealth}]
 \begin{scope}[scale=1.6]
  \fill[blue!7] (0,0) circle (1);
  \draw (0,0) circle (1);
  \draw[->] (-1.3,0)--(1.4,0) node[right] {$x$};
  \draw[->] (0,-1.3)--(0,1.4) node[above] {$y$};
  \draw[very thick,blue] (0.5,-0.866)--(0.5,0.866);
  \node[below right] at (0.5,0) {$t$};
  \node[blue,right] at (0.52,0.55) {chord at $x=t$};
 \end{scope}
 \begin{scope}[xshift=6.2cm,x=1.6cm,y=1cm]
  \draw[->] (-1.3,0)--(1.4,0) node[right] {$t$};
  \draw[->] (0,-0.2)--(0,2.5) node[above] {$G(t)$};
  \draw[thick,blue,domain=0:180,samples=90,smooth,variable=\a] plot ({cos(\a)},{2*sin(\a)});
  \node[below] at (-1,0) {$-1$};
  \node[below] at (1,0) {$1$};
  \draw (-0.05,2)--(0.05,2);
  \node[above left] at (0,2) {$2$};
 \end{scope}
\end{tikzpicture}
\caption{The chord integrals of $u\equiv1$ on the unit disk. The length $G(t)=2\sqrt{1-t^2}$ of the chord at $x=t$ is a strictly concave function of $t$, because the endpoints of the chord move with $t$ and their contributions survive differentiation.}
\label{fig:chords}
\end{figure}
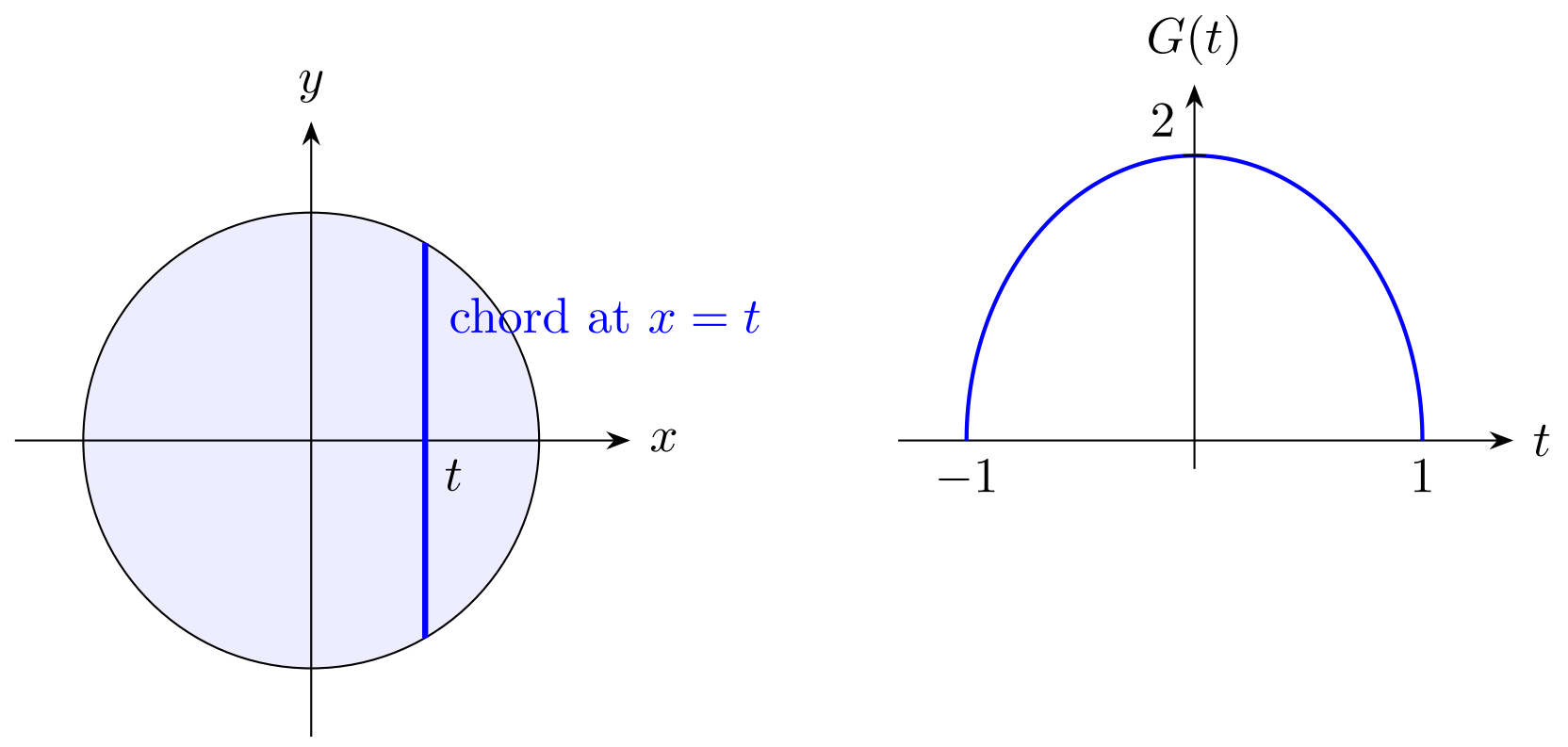

Let us now see what changes when we integrate over the chords of a bounded domain. Let $B=\{(x,y)\in\mathbb{R}^2\mid x^2+y^2<1\}$ and take $u\equiv1$ on $B$, a harmonic function on $B$. Its integral over the vertical chord at $x=t$ is $G(t)=\int_{-\sqrt{1-t^2}}^{\sqrt{1-t^2}}dy=2\sqrt{1-t^2}$ for $|t|<1$, and yet $G''(t)=-\frac{2}{(1-t^2)^{\frac{3}{2}}}<0$, so the chord integral is strictly concave, as Figure~\ref{fig:chords} shows. The endpoints of the chord now depend on $t$, and their contributions remain when we differentiate. Extending $u$ by zero outside the disk does not repair the argument, because the extension is not subharmonic on the whole plane. To see this directly, take a nonnegative smooth radial test function $\psi(x)=\eta(|x|)$ with compact support and $\eta'(1)<0$. By the divergence theorem, $\int_B\Delta\psi(x)dx=2\pi\eta'(1)<0$, so the extension violates \eqref{eq_subharmonic}.

Zero boundary values do not help either. The function $u(x,y)=x^2+y^2-1$ is subharmonic on $B$ and vanishes on its boundary, but its chord integrals are not convex. With $\theta=(1,0)$ and $h=\sqrt{1-t^2}$ we have $t^2-1=-h^2$, so
$$
Ru(\theta,t)=\int_{-h}^h(y^2-h^2)dy=-\frac{4}{3}h^3=-\frac{4}{3}(1-t^2)^{\frac{3}{2}},
$$
and differentiating twice, we find that $\frac{d^2}{dt^2}Ru(\theta,t)=\frac{4(1-2t^2)}{\sqrt{1-t^2}}$, which is negative when $\frac{1}{2}<t^2<1$. The Radon transform is therefore not convex on $(-1,1)$.

For a smooth function $u$ on the closure of a smooth bounded domain $D$, the boundary contribution can be written explicitly. Let $1_D$ be the indicator function of $D$, let $\nu$ be the outward unit normal, and let $\partial_\nu$ denote differentiation in the direction $\nu$. For a distribution $T$ and a test function $\psi$ we write $\langle T,\psi\rangle$ for the action of $T$ on $\psi$, and Green's formula reads
\begin{equation}\label{eq_boundary}
\langle\Delta(u1_D),\psi\rangle=\int_D(\Delta u)\psi dx+\int_{\partial D}(u\partial_\nu\psi-(\partial_\nu u)\psi)dS.
\end{equation}
It follows from two integrations by parts. The first one shows that
$$
\langle\Delta(u1_D),\psi\rangle=\int_Du\Delta\psi dx=\int_{\partial D}u\partial_\nu\psi dS-\sum_{j=1}^d\int_D(\partial_ju)(\partial_j\psi)dx,
$$
and integrating the last terms by parts once more produces $\int_D(\Delta u)\psi dx$ and $-\int_{\partial D}(\partial_\nu u)\psi dS$, which proves \eqref{eq_boundary}. Neither boundary term can be discarded without an additional assumption.

There is also a limitation in the opposite direction. If $u$ is subharmonic on all of $\mathbb R^d$ and belongs to $L^1(\mathbb R^d)$, then $u=0$ almost everywhere, so on the whole space integrability leaves no room for a nonzero subharmonic function at all. Indeed, for $u\in L^1(\mathbb R^d)$ the hypothesis \eqref{eq_slab} holds in every direction, so by Theorem~\ref{thm_weak} every directional Radon transform agrees almost everywhere with a convex function on $\mathbb R$, which is integrable by Fubini's theorem, and such a function must vanish. To see this, suppose that a convex function $G$ has a positive secant slope. Then all secant slopes to its right are at least as large, which forces linear growth to the right, and similarly a negative secant slope forces linear growth to the left. Both possibilities contradict integrability, so every secant slope is zero, $G$ is constant, and integrability forces this constant to be zero. The Fourier slice theorem and the uniqueness of the Fourier transform now show that $u=0$. On a slab there is no such obstruction, since, for example, the restriction of $e^{\sqrt2x-y^2}$ to $a<x<b$ is integrable on that slab and has a nonzero convex Radon transform in the direction $(1,0)$.

\section{One-dimensional Fourier transform}\label{odres}

Before applying the Fourier slice theorem, we collect the one-dimensional results that we shall need. We begin with the balance operator, then derive the Fourier representation, and finally prove the remainder estimate for convex functions, which is the part that enters the multidimensional argument. We refer to \cite{Ltheta,L2019,T0,TB} for the background.

\subsection{Balance integral operator}\label{balance}

For $x>0$ and a function $\varphi$ on $(0,\infty)$, the balance operator is defined by
$$
B_\varphi g(x)=\frac{1}{x^2}\int_0^\infty g\left(\frac{t}{x}\right)\varphi(t)dt.
$$
A change of variables shows that
$$
xB_{\sin(\cdot)}g(x)=\int_0^\infty g(t)\sin(xt)dt
$$
and
$$
xB_{\cos(\cdot)}g(x)=\int_0^\infty g(t)\cos(xt)dt
$$
whenever the integrals exist, so the Fourier transform on the whole line can be treated by separating the positive and the negative half axis. We shall use the following elementary estimate.

\begin{lemma}\label{mp1}
Let $g\in L^1(\mathbb R_+)$, where $\mathbb R_+=(0,\infty)$. If $\int_0^\infty\frac{|\varphi(t)|}{t}dt<\infty$, then $B_\varphi g\in L^1(\mathbb R_+)$ and
$$
\int_0^\infty|B_\varphi g(x)|dx\leq\int_0^\infty|g(s)|ds\int_0^\infty\frac{|\varphi(t)|}{t}dt.
$$
\end{lemma}

\begin{proof}
By Tonelli's theorem and the substitution $s=\frac{t}{x}$,
$$
\int_0^\infty\int_0^\infty\frac{1}{x^2}\left|g\left(\frac{t}{x}\right)\right||\varphi(t)|dt dx=
\int_0^\infty\frac{|\varphi(t)|}{t}dt\int_0^\infty|g(s)|ds,
$$
and the assertion follows.
\end{proof}

On a finite interval $[a,b]$ we use the same notation for the operator
\begin{equation}\label{eq_balance_finite}
B_\varphi g(x)=\frac{1}{|x|^2}\int_{a|x|}^{b|x|}g\left(\frac{s}{|x|}\right)\varphi(s,x)ds,
\end{equation}
defined for $x\ne0$ and for a kernel $\varphi(s,x)$ that will be specified below. The absolute values in the limits and in the argument of $g$ allow us to use this definition for both signs of $x$.

\subsection{General representation of the one-dimensional Fourier transform}

We keep the sign convention of \eqref{gft} and write
\begin{equation}\label{nftab}
\widehat f(x)=\int_a^bf(t)e^{-ixt}dt
\end{equation}
for a function $f$ on a finite interval $[a,b]$. Here and below $W^{1,1}(a,b)$ is the space of integrable functions $f$ on $(a,b)$ with an integrable weak derivative, that is, an integrable function $g$ such that $\int_a^bf(t)\psi'(t)dt=-\int_a^bg(t)\psi(t)dt$ for every $\psi\in C_c^\infty(a,b)$, and we write $f'=g$. Each such function has an absolutely continuous representative on $[a,b]$, and point values always refer to this representative.

\begin{figure}[htbp]
\centering
\begin{tikzpicture}[x=1cm,y=1cm,font=\small,>={Stealth}]
 \fill[blue!12] (0,0) rectangle (1.6,0.45);
 \fill[orange!20] (1.6,0) rectangle (8.4,0.45);
 \fill[blue!12] (8.4,0) rectangle (10,0.45);
 \draw[thick] (0,0)--(10,0);
 \foreach \x/\l in {0/a,1.6/a+h,8.4/b-h,10/b} {\draw (\x,0.12)--(\x,-0.12) node[below] {$\l$};}
 \node[above] at (0.8,0.5) {$e^{-ixt}-e^{-iax}$};
 \node[above] at (5,0.5) {$e^{-ixt}$};
 \node[above] at (9.2,0.5) {$e^{-ixt}-e^{-ibx}$};
 \draw[<->] (0,-0.95)--(1.6,-0.95);
 \node[below] at (0.8,-1.0) {$h=\frac{d}{|x|}$};
 \draw[<->] (8.4,-0.95)--(10,-0.95);
 \node[below] at (9.2,-1.0) {$h$};
 \node[below,align=center,font=\footnotesize] at (5,-0.55) {monotone $f'$ when $f$ is convex,\\ second integration by parts};
\end{tikzpicture}
\caption{The decomposition in the proof of Theorem~\ref{odftc1}. The endpoint values are moved a distance $h=\frac{d}{|x|}$ into $[a,b]$, and above each of the three intervals we show the factor that multiplies $f'$ in the remainder. On the two short intervals this factor is a difference of exponentials of size at most $|x|h=d$, and on the middle interval, when $f$ is convex, $f'$ is monotone and a second integration by parts applies.}
\label{fig:three}
\end{figure}
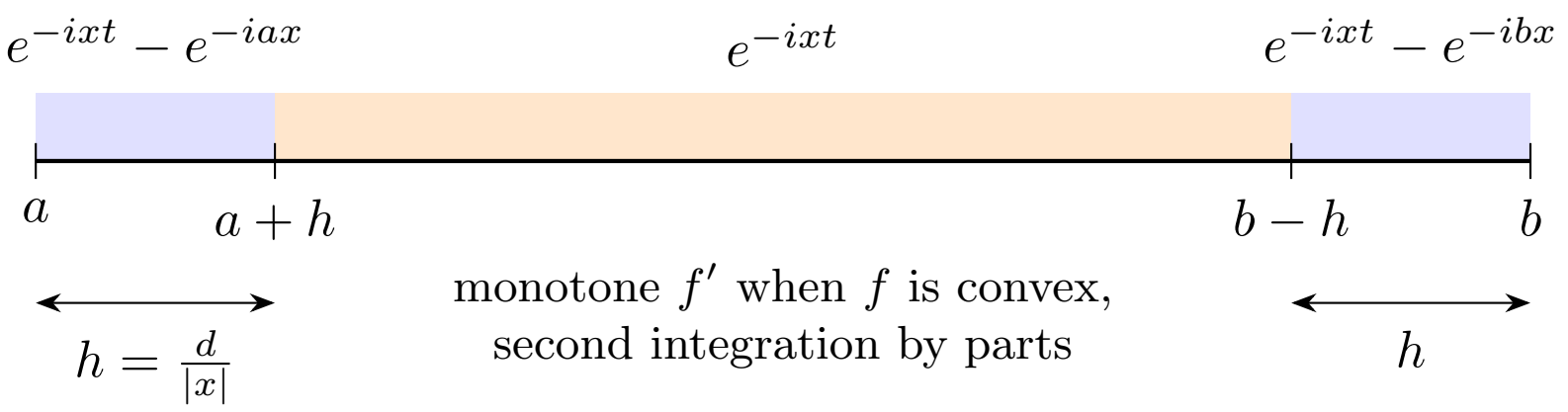

Let us first see what integration by parts does. It produces the two endpoint values of $f$ divided by $x$, together with an oscillatory integral of $f'$, also divided by $x$. Taking absolute values in this last integral would leave a bound of order $\frac{1}{|x|}$, which is not integrable at infinity, so a second idea is needed. We move the endpoint values a distance $\frac{d}{|x|}$ into the interval, as in Figure~\ref{fig:three}. The errors near the endpoints then contain differences of exponentials, which are small, while the middle integral can take advantage of monotonicity when $f$ is convex.

\begin{theorem}\label{odftc1}
Let $-\infty<a<b<\infty$, let $f\in W^{1,1}(a,b)$, and put $d=\min\{\frac{b-a}{2},\frac{\pi}{2}\}$. Then, for $|x|>2$,
\begin{equation}\label{ftgmain1}
\widehat f(x)=-\frac{i}{x}f\left(a+\frac{d}{|x|}\right)e^{-iax}+\frac{i}{x}f\left(b-\frac{d}{|x|}\right)e^{-ibx}+B_\varphi(f')(x),
\end{equation}
where $B_\varphi$ is defined by \eqref{eq_balance_finite} with the kernel
$$
\varphi(s,x)=-i\operatorname{sign}x(e^{-is\operatorname{sign}x}-e^{-iax}1_{[a|x|,a|x|+d)}(s)-e^{-ibx}1_{(b|x|-d,b|x|]}(s))1_{[a|x|,b|x|]}(s),
$$
and $1_I$ denotes the indicator function of an interval $I$.
\end{theorem}

\begin{proof}
Put $h=\frac{d}{|x|}$. Since $|x|>2$ and $d\leq\frac{b-a}{2}$, we have $h<\frac{b-a}{4}$, so the points $a<a+h<b-h<b$ are ordered, and the two endpoint intervals in the definition of the kernel are disjoint. Integration by parts shows that
$$
\widehat f(x)=-\frac{i}{x}f(a)e^{-iax}+\frac{i}{x}f(b)e^{-ibx}-\frac{i}{x}\int_a^bf'(t)e^{-ixt}dt,
$$
and we insert the identities $f(a)=f(a+h)-\int_a^{a+h}f'(t)dt$ and $f(b)=f(b-h)+\int_{b-h}^bf'(t)dt$. After collecting the terms, the remainder is $-\frac{i}{x}$ times
\begin{equation}\label{eq_three_terms}
\int_a^{a+h}f'(t)(e^{-ixt}-e^{-iax})dt+\int_{a+h}^{b-h}f'(t)e^{-ixt}dt+\int_{b-h}^bf'(t)(e^{-ixt}-e^{-ibx})dt,
\end{equation}
and the substitution $s=|x|t$ identifies this expression, multiplied by $-\frac{i}{x}$, with $B_\varphi(f')(x)$.
\end{proof}

The representation holds without any convexity, but the reason for separating the three intervals becomes apparent when $f$ is convex. Near each endpoint the difference of exponentials is small, and in the middle, convexity makes the derivative monotone, which allows a second integration by parts. The resulting estimate is associated with the work of Trigub \cite{T0,TB}, with earlier results for Fourier series due to Shilov \cite{Sh}, and we include its proof in the form needed here.

\begin{corollary}\label{cor_convex_remainder}
If $f\in W^{1,1}(a,b)$ is convex, then the remainder in \eqref{ftgmain1} satisfies
\begin{equation}\label{rftab}
\int_{|x|>2}|B_\varphi(f')(x)|dx
\leq C\left(d+\frac{1}{d}\right)\int_a^b|f'(t)|dt.
\end{equation}
\end{corollary}

\begin{proof}
Let $g=f'_+$ be the right derivative of $f$, which is nondecreasing and right-continuous on $(a,b)$ and agrees with $f'$ almost everywhere. We estimate the three terms in \eqref{eq_three_terms}, each multiplied by $-\frac{i}{x}$, and write $r=|x|>2$.

Since $|e^{-ixt}-e^{-iax}|\leq r(t-a)$, the first term is at most $\int_a^{a+\frac{d}{r}}(t-a)|g(t)|dt$ in absolute value. For a fixed $v=t-a$, the variable $r$ ranges over $2<r<\frac{d}{v}$, so Tonelli's theorem bounds the integral of this quantity over $r>2$ by
$$
\int_0^{\frac{d}{2}}v|g(a+v)|\left(\frac{d}{v}-2\right)dv\leq d\int_0^{\frac{d}{2}}|g(a+v)|dv.
$$
The same argument at $b$ bounds the integral of the third term over $r>2$ by $d\int_0^{\frac{d}{2}}|g(b-v)|dv$, and these estimates hold for both signs of $x$.

For the middle term, put $\alpha=a+\frac{d}{r}$ and $\beta=b-\frac{d}{r}$. The monotone function $g$ has finite variation on $[\alpha,\beta]$, and its derivative is the positive measure $dg$ determined by $dg((c,e])=g(e)-g(c)$ for $\alpha\leq c<e\leq\beta$. For smooth $g$ this is simply $g'(t)dt$. For a general monotone $g$, the identity $g(t)=g(\alpha)+dg((\alpha,t])$ and Fubini's theorem lead to the integration by parts formula
$$
\int_\alpha^\beta g(t)e^{-ixt}dt=\frac{i}{x}(g(\beta)e^{-ix\beta}-g(\alpha)e^{-ix\alpha})-\frac{i}{x}\int_{(\alpha,\beta]}e^{-ixt}dg(t).
$$
The total mass of $dg$ on $(\alpha,\beta]$ is $g(\beta)-g(\alpha)$, and taking absolute values we obtain
$$
\left|\int_\alpha^\beta g(t)e^{-ixt}dt\right|\leq\frac{|g(\alpha)|+|g(\beta)|+g(\beta)-g(\alpha)}{r}\leq\frac{2}{r}(|g(\alpha)|+|g(\beta)|).
$$
This argument only uses values of $g$ strictly inside $(a,b)$, so no boundedness of $g$ at $a$ or $b$ is needed. Including the factor $\frac{1}{r}$, the integral of the middle term over either half axis is therefore at most
$$
\int_2^\infty\frac{2}{r^2}\left(\left|g\left(a+\frac{d}{r}\right)\right|+
\left|g\left(b-\frac{d}{r}\right)\right|\right)dr=\frac{2}{d}\int_0^{\frac{d}{2}}(|g(a+v)|+|g(b-v)|)dv,
$$
where we substituted $v=\frac{d}{r}$. The intervals $(a,a+\frac{d}{2})$ and $(b-\frac{d}{2},b)$ are disjoint, so summing the estimates for both signs of $x$ proves \eqref{rftab} with $2d+\frac{4}{d}$ in place of $C(d+\frac{1}{d})$.
\end{proof}

A finite convex function on $(a,b)$ is locally absolutely continuous, but its derivative need not be integrable on the whole interval. For example, $f(t)=\frac{1}{\sqrt t}$ is convex and belongs to $L^1(0,1)$, whereas
$$
\int_0^1|f'(t)|dt=\frac{1}{2}\int_0^1\frac{dt}{t^{\frac{3}{2}}}=\infty.
$$
The assumption $f\in W^{1,1}(a,b)$ in Corollary~\ref{cor_convex_remainder} therefore includes control at the endpoints. It holds, in particular, for a convex function that is continuous on $[a,b]$.

\section{General representation of the multivariate Fourier transform}\label{gerai}

We call an open convex domain $D\subset\mathbb R^d$ containing the origin admissible. If $f\in L^1(D)$, its Fourier and Radon transforms are those of its extension by zero outside $D$. For a fixed direction $\theta\in\mathbb S^{d-1}$ put $a_\theta=\inf_{u\in D}\langle u,\theta\rangle$ and $b_\theta=\sup_{u\in D}\langle u,\theta\rangle$, and whenever these numbers are finite, set $d_\theta=\min\{\frac{b_\theta-a_\theta}{2},\frac{\pi}{2}\}$. The support of $Rf(\theta,\cdot)$ lies in $[a_\theta,b_\theta]$. For $|t|>2$ we write
$$
E_\theta(t)=-\frac{i}{t}Rf\left(\theta,a_\theta+\frac{d_\theta}{|t|}\right)e^{-ia_\theta t}+\frac{i}{t}Rf\left(\theta,b_\theta-\frac{d_\theta}{|t|}\right)e^{-ib_\theta t}
$$
for the sum of the two endpoint terms that will appear in the directional formulas.

\begin{theorem}\label{odftm}
Let $f\in L^1(D)$, where $D$ is admissible, and fix $\theta\in\mathbb S^{d-1}$ such that $a_\theta$ and $b_\theta$ are finite. Suppose that $Rf(\theta,\cdot)\in W^{1,1}(a_\theta,b_\theta)$, and let point values of $Rf(\theta,\cdot)$ refer to its absolutely continuous representative. Then, for $|t|>2$,
\begin{equation}\label{ftgmain}
\widehat f(t\theta)=E_\theta(t)+B_{\varphi_\theta}\left(\frac{\partial Rf}{\partial s}(\theta,\cdot)\right)(t),
\end{equation}
where the balance operator is defined by \eqref{eq_balance_finite} with $a=a_\theta$ and $b=b_\theta$, and $\varphi_\theta$ is the kernel of Theorem~\ref{odftc1} with $a$, $b$, $d$ and $x$ replaced by $a_\theta$, $b_\theta$, $d_\theta$ and $t$.
\end{theorem}

\begin{proof}
By the Fourier slice theorem \eqref{t5.9}, $\widehat f(t\theta)$ is the one-dimensional Fourier transform of $Rf(\theta,\cdot)$, which is supported in $[a_\theta,b_\theta]$, and Theorem~\ref{odftc1} on $(a_\theta,b_\theta)$, with $x$ replaced by $t$, proves \eqref{ftgmain}.
\end{proof}

Theorem~\ref{odftm} is a representation for a general integrable function whose Radon transform has an integrable first derivative in the chosen direction. To control the remainder we need an additional property of this Radon transform, and this is where the convexity of Sections~\ref{sec_smooth} and~\ref{sec_weak} enters.

\subsection{Subharmonic functions}

\begin{theorem}\label{corsbh}
Fix $\theta\in\mathbb S^{d-1}$ and finite numbers $a_\theta<0<b_\theta$, and let $f\in L^1(D)$ be subharmonic on the slab $D=\{u\in\mathbb R^d\mid a_\theta<\langle u,\theta\rangle<b_\theta\}$. This slab is an admissible domain, the numbers $a_\theta$ and $b_\theta$ of the notation introduced above are the given ones, and $d_\theta=\min\{\frac{b_\theta-a_\theta}{2},\frac{\pi}{2}\}$. Let $Rf(\theta,\cdot)$ denote the convex representative supplied by Corollary~\ref{cor_slab}, and suppose that the limits $A_\theta=\lim_{s\downarrow a_\theta}Rf(\theta,s)$ and $B_\theta=\lim_{s\uparrow b_\theta}Rf(\theta,s)$ are finite. Then the derivative of $Rf(\theta,\cdot)$ is integrable on $(a_\theta,b_\theta)$, and, for $|t|>2$,
\begin{equation}\label{corftgmain}
\widehat f(t\theta)=E_\theta(t)+F(\theta,t),
\end{equation}
where $E_\theta(t)$ is computed with this representative and
\begin{equation}\label{estrems}
\int_{|t|>2}|F(\theta,t)|dt\leq C\left(d_\theta+\frac{1}{d_\theta}\right)\int_{a_\theta}^{b_\theta}\left|\frac{\partial}{\partial s}Rf(\theta,s)\right|ds.
\end{equation}
\end{theorem}

\begin{proof}
Since $f\in L^1(D)$, the integrability condition of Corollary~\ref{cor_slab} holds, and $Rf(\theta,\cdot)$ has a finite convex representative on $(a_\theta,b_\theta)$. Write $G(s)=Rf(\theta,s)$. The finite endpoint limits extend $G$ to a continuous convex function on $[a_\theta,b_\theta]$, with $G(a_\theta)=A_\theta$ and $G(b_\theta)=B_\theta$. Choose a point $c_\theta$ where this extension attains its minimum $m_\theta$. The secant slope inequalities show that $G$ is nonincreasing to the left of $c_\theta$ and nondecreasing to its right.

On each compact subinterval of $(a_\theta,b_\theta)$, the function $G$ is absolutely continuous. Integrating $-G'$ to the left of $c_\theta$ and $G'$ to its right, and then letting the endpoints of these subintervals approach $a_\theta$ and $b_\theta$, gives
\begin{equation}\label{eq_slab_variation}
\int_{a_\theta}^{b_\theta}|G'(s)|ds
=A_\theta+B_\theta-2m_\theta<\infty.
\end{equation}
If the minimum occurs at an endpoint, the same calculation uses only the other interval, and in either case the derivative integral is finite because of convexity and the finite endpoint limits.

Being convex, the representative is locally absolutely continuous, so for every $c\in(a_\theta,b_\theta)$
$$
Rf(\theta,s)=Rf(\theta,c)+\int_c^s\frac{\partial}{\partial v}Rf(\theta,v)dv.
$$
The derivative is integrable on the whole interval, so the right side has finite limits at both endpoints and defines an absolutely continuous extension to $[a_\theta,b_\theta]$. In particular $Rf(\theta,\cdot)\in W^{1,1}(a_\theta,b_\theta)$, Theorem~\ref{odftm} identifies $F(\theta,t)$ with the balance remainder, and \eqref{rftab} proves \eqref{estrems}.
\end{proof}

\begin{remark}
No continuity assumption on $f$ is needed, since the direction $\theta$ is fixed and the convex representative agrees with the Radon integral for almost every $s$. Nothing is claimed about other directions.
\end{remark}

The finite endpoint limits are equivalent to integrability of the derivative for a finite convex function on a bounded interval. One implication was proved in \eqref{eq_slab_variation}. For the other, the identity in the preceding proof gives finite endpoint limits whenever the derivative is integrable.

Subharmonicity and integrability on the slab do not by themselves give these endpoint limits. In the plane, let $D=\{(s,y)\in\mathbb R^2\mid -1<s<1\}$ and set
$$
f(s,y)=\frac{e^{-\frac{y^2}{16}}}{\sqrt{1-s}}.
$$
This function belongs to $L^1(D)$, and
$$
\Delta f(s,y)=\left(\frac{3}{4(1-s)^2}+\frac{y^2}{64}-\frac{1}{8}\right)f(s,y)
\geq\frac{1}{16}f(s,y)>0.
$$
For $\theta=(1,0)$ its Radon transform is $Rf(\theta,s)=\frac{4\sqrt{\pi}}{\sqrt{1-s}}$. This is convex and integrable on $(-1,1)$, but
$$
\int_{-1}^1\left|\frac{\partial}{\partial s}Rf(\theta,s)\right|ds
=2\sqrt{\pi}\int_{-1}^1\frac{ds}{(1-s)^{\frac{3}{2}}}=\infty.
$$
The divergence occurs at the endpoint $s=1$.

A useful sufficient condition for the endpoint limits in Theorem~\ref{corsbh} is $f\in W^{1,1}(D)$, where $W^{1,1}(D)$ is the space of integrable functions on $D$ whose first weak derivatives are integrable. The weak derivatives $\partial_jf$ are defined by the identities $\int_Df\partial_j\psi=-\int_D(\partial_jf)\psi$ for every $\psi\in C_c^\infty(D)$ and every coordinate $j$, and $\triangledown f$ denotes the vector $(\partial_1f,\ldots,\partial_df)$.

We claim that, for $f\in W^{1,1}(D)$,
\begin{equation}\label{eq_slab_derivative}
\frac{\partial}{\partial s}Rf(\theta,s)=\int_{\theta^\perp}\langle\theta,\triangledown f(s\theta+y)\rangle dy
\end{equation}
for almost every $s\in(a_\theta,b_\theta)$. To verify this, take $\varphi\in C_c^\infty(a_\theta,b_\theta)$ and use the test function $\varphi(s)\chi_R(y)$, with the cutoff $\chi_R$ of Section~\ref{sec_smooth}. Since $\chi_R$ does not depend on $s$, the derivative of this test function in the direction $\theta$ is $\varphi'(s)\chi_R(y)$, and the weak derivative identity reads
$$
\int_Df(s\theta+y)\varphi'(s)\chi_R(y)dy ds=-\int_D\langle\theta,\triangledown f(s\theta+y)\rangle\varphi(s)\chi_R(y)dy ds.
$$
Both $f$ and its gradient are integrable, so dominated convergence allows us to let $R\to\infty$, and Fubini's theorem identifies the left side with $\int_{a_\theta}^{b_\theta}Rf(\theta,s)\varphi'(s)ds$ and the right side with minus the integral of the proposed derivative against $\varphi$. This proves \eqref{eq_slab_derivative}, and taking absolute values and integrating once more, we arrive at
$$
\int_{a_\theta}^{b_\theta}\left|\frac{\partial}{\partial s}Rf(\theta,s)\right|ds\leq\int_D|\langle\theta,\triangledown f(u)\rangle|du\leq\int_D|\triangledown f(u)|du.
$$
Together with $Rf(\theta,\cdot)\in L^1(a_\theta,b_\theta)$, this puts the Radon transform in $W^{1,1}(a_\theta,b_\theta)$ and gives finite endpoint limits. The directional remainder is thus controlled by a quantity involving $f$ itself, which is what we were after.

\subsection{Integrability of the Fourier transform}

Averaging a directional estimate over $\theta$ requires attention to the measure in frequency space. In polar coordinates,
\begin{equation}\label{eq_polar_weight}
\int_{|x|>2}\frac{|\widehat f(x)|}{|x|^{d-1}}dx=\int_{\mathbb S^{d-1}}\int_2^\infty|\widehat f(t\theta)|dt d\theta,
\end{equation}
where $d\theta$ is surface measure on $\mathbb S^{d-1}$, so the unweighted integral in the radial variable corresponds to a weighted integral in $\mathbb R^d$.

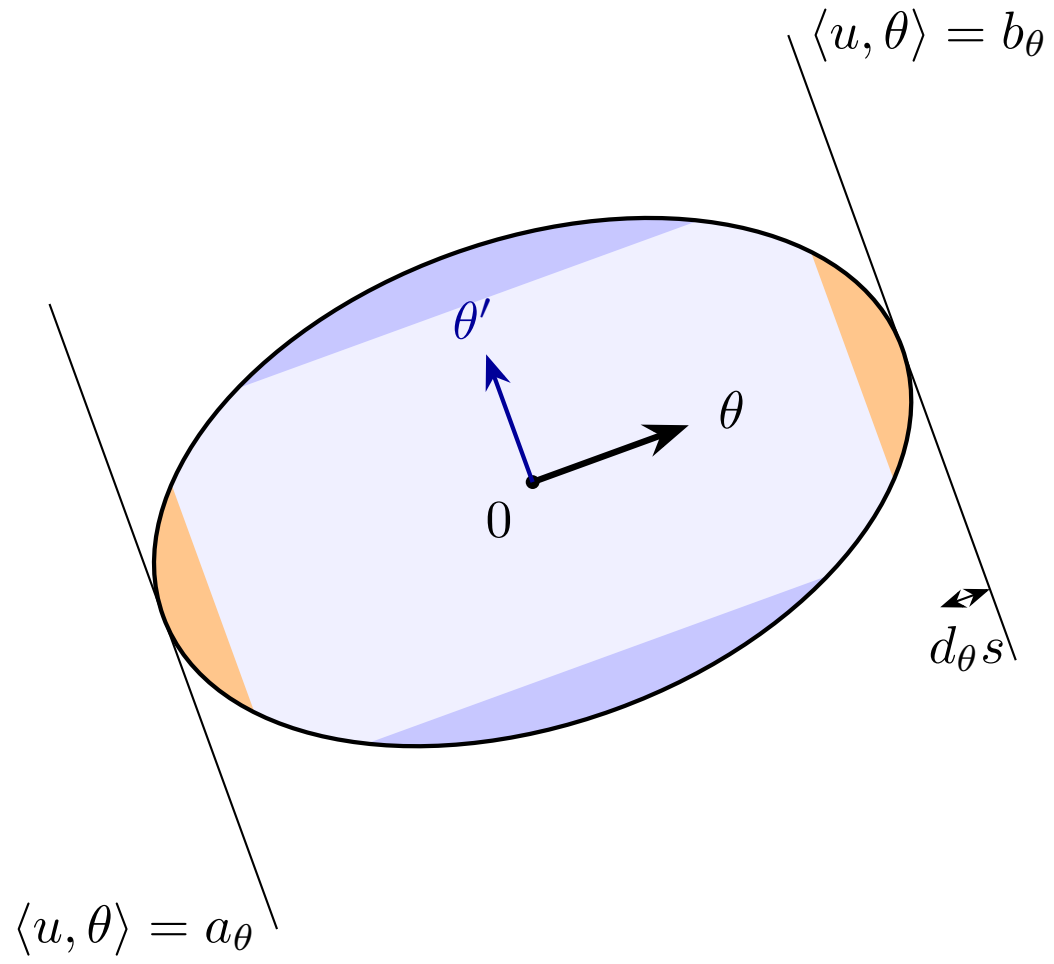
\begin{figure}[htbp]
\centering
\begin{tikzpicture}[font=\small,>={Stealth}]
 \begin{scope}
  \clip[rotate=20] (0,0) ellipse (2.6 and 1.6);
  \fill[blue!6,rotate=20] (0,0) ellipse (2.6 and 1.6);
  \fill[blue!22,rotate=20] (-3.5,1.25) rectangle (3.5,1.6);
  \fill[blue!22,rotate=20] (-3.5,-1.6) rectangle (3.5,-1.25);
  \fill[orange!45,rotate=20] (2.25,-2.5) rectangle (2.6,2.5);
  \fill[orange!45,rotate=20] (-2.6,-2.5) rectangle (-2.25,2.5);
 \end{scope}
 \draw[thick,rotate=20] (0,0) ellipse (2.6 and 1.6);
 \fill (0,0) circle (1.3pt);
 \node[below left] at (0,0) {$0$};
 \draw[->,very thick] (0,0)--(20:1.1);
 \node at (20:1.4) {$\theta$};
 \draw[->,thick,blue!60!black] (0,0)--(110:0.9);
 \node[blue!60!black] at (110:1.15) {$\theta'$};
 \draw ($(20:2.6)+(110:2.2)$)--($(20:2.6)+(-70:2.2)$);
 \draw ($(200:2.6)+(110:2.2)$)--($(200:2.6)+(-70:2.2)$);
 \node[right] at ($(20:2.6)+(110:2.2)$) {$\langle u,\theta\rangle=b_\theta$};
 \node[left] at ($(200:2.6)+(-70:2.2)$) {$\langle u,\theta\rangle=a_\theta$};
 \draw[<->] ($(20:2.25)+(-70:1.7)$)--($(20:2.6)+(-70:1.7)$);
 \node[below] at ($(20:2.42)+(-70:1.75)$) {$d_\theta s$};
\end{tikzpicture}
\caption{The boundary layer of a bounded convex domain. For the direction $\theta$ the supporting hyperplanes $\langle u,\theta\rangle=a_\theta$ and $\langle u,\theta\rangle=b_\theta$ cut off the two darker strips of width $d_\theta s$ inside $D$, and the lighter strips are the contribution of a second direction $\theta'$. The set $D_s$ is the union of such strips over all directions, a layer along $\partial D$ whose thickness lies between $d_*s$ and $d^*s$.}
\label{fig:layer}
\end{figure}

Let $D$ now be a bounded admissible domain. For $0<s<1$ we define its boundary layer $D_s$ as the set of points $u\in D$ for which there is a direction $\theta\in\mathbb S^{d-1}$ with
\begin{equation}\label{eq_layer}
0<\langle u,\theta\rangle-a_\theta<d_\theta s\text{ or }0<b_\theta-\langle u,\theta\rangle<d_\theta s.
\end{equation}
Figure~\ref{fig:layer} shows the strips that a single direction contributes to $D_s$. The numbers $d_\theta$ are bounded above and below by positive constants depending on $D$. Indeed, if the ball of radius $r_0>0$ about the origin lies in $D$, then $b_\theta-a_\theta\geq2r_0$, and therefore $\min\{r_0,\frac{\pi}{2}\}\leq d_\theta\leq\frac{\pi}{2}$. Boundedness of $D$ also makes $a_\theta$ and $b_\theta$ continuous in $\theta$, since $|u|\leq R_0$ on $D$ implies $|b_\theta-b_\omega|\leq R_0|\theta-\omega|$, and similarly for $a_\theta$. Hence $d_\theta$ is continuous on the sphere and attains its minimum $d_*$ and its maximum $d^*$. The set $D_s$ is open, and it lies near the boundary, in the sense that
$$
\{u\in D\mid\operatorname{dist}(u,\partial D)<d_*s\}\subset D_s\subset\{u\in D\mid\operatorname{dist}(u,\partial D)<d^*s\}.
$$
To justify these inclusions, fix $u\in D$. A supporting hyperplane of $D$ is a hyperplane that bounds a half space containing $D$ and touches the closure of $D$, and the distances from $u$ to the two supporting hyperplanes normal to $\theta$ are $\langle u,\theta\rangle-a_\theta$ and $b_\theta-\langle u,\theta\rangle$. The ball centered at $u$ with radius $\operatorname{dist}(u,\partial D)$ lies in $D$, so all these distances are at least this radius. Conversely, the supporting hyperplane through a nearest boundary point is at distance at most $\operatorname{dist}(u,\partial D)$ from $u$. Hence the smallest distance from $u$ to a supporting hyperplane equals $\operatorname{dist}(u,\partial D)$, and the two inclusions follow from $d_*\leq d_\theta\leq d^*$.

On a bounded domain the examples of Section~\ref{sec_examples} prevent us from deducing convexity of $Rf$ from subharmonicity alone, so we state this assumption explicitly. We also impose a condition on the extension by zero, so that differentiation does not produce a measure on the boundary.

\begin{theorem}\label{avft}
Let $D$ be a bounded admissible domain, and let $f\in W^{1,1}(D)$. Suppose that its extension by zero belongs to $W^{1,1}(\mathbb R^d)$ and that, for almost every $\theta\in\mathbb S^{d-1}$, $Rf(\theta,\cdot)$ has a convex representative on $(a_\theta,b_\theta)$. Then
\begin{equation}\label{estintft}
\int_{|x|>2}\frac{|\widehat f(x)|}{|x|^{d-1}}dx
\leq C_D\int_0^1\frac{1}{t^2}\int_{D_t}|f(u)|du dt
+C_D\int_D|\triangledown f(u)|du.
\end{equation}
In particular, the weighted Fourier integral is finite whenever the boundary layer integral on the right is finite.
\end{theorem}

\begin{proof}
As stipulated at the beginning of this section, $Rf$ denotes the Radon transform of the extension of $f$ by zero. The weak gradient of this extension agrees with $\triangledown f$ inside $D$ and vanishes outside $\overline D$, and since the boundary of a bounded convex domain has Lebesgue measure zero, the integral of its absolute value is exactly $\int_D|\triangledown f|$. The test function argument that proved \eqref{eq_slab_derivative}, applied to the extension, shows that
$$
\frac{\partial}{\partial s}Rf(\theta,s)=R(\langle\theta,\triangledown f\rangle1_D)(\theta,s)
$$
in the sense of distributions and almost everywhere, and in particular
\begin{equation}\label{eq_radon_gradient}
\int_{a_\theta}^{b_\theta}\left|\frac{\partial}{\partial s}Rf(\theta,s)\right|ds\leq\int_D|\triangledown f(u)|du.
\end{equation}
For smooth functions on the whole space this is the first derivative case of the usual differentiation identity for the Radon transform \cite{Helgason,Natterer}, and the weak formulation is needed here because the extension need not be smooth.

By Fubini's theorem we also have $Rf(\theta,\cdot)\in L^1(\mathbb R)$, and together with the derivative identity this puts $Rf(\theta,\cdot)$ in $W^{1,1}(\mathbb R)$. Its continuous representative agrees everywhere on $(a_\theta,b_\theta)$ with the convex representative whenever the latter exists, so Theorem~\ref{odftm} and Corollary~\ref{cor_convex_remainder} apply in every direction in which $Rf(\theta,\cdot)$ is convex, with point values taken from the convex representative. Changing $Rf(\theta,\cdot)$ on a null set does not affect any of the integrals below, since the changes of variables are smooth with nonzero derivatives on the intervals involved. To make the choice of representatives precise when we average over directions, we take
$$
Rf(\theta,s)=\lim_{n\to\infty}\int_Df(u)\eta_{\frac{1}{n}}(s-\langle u,\theta\rangle)du,
$$
where $\eta_{\frac{1}{n}}$ is the approximate identity of Section~\ref{sec_weak}. The integrals on the right are jointly continuous in $(\theta,s)$ by dominated convergence, and for each fixed $\theta$ their limit is the continuous representative. The endpoint terms $E_\theta(t)$ are therefore measurable in the direction and the frequency, and so is $F(\theta,t)=\widehat f(t\theta)-E_\theta(t)$. The exceptional directions, in which convexity is not assumed, form a set of zero surface measure and do not affect \eqref{eq_polar_weight}.

We first estimate the term involving $a_\theta$. The change of variables $h=\frac{d_\theta}{t}$ shows that
$$
\int_2^\infty\frac{1}{t}\left|Rf\left(\theta,a_\theta+\frac{d_\theta}{t}\right)\right|dt=
\int_0^{\frac{d_\theta}{2}}\frac{|Rf(\theta,a_\theta+h)|}{h}dh,
$$
and integrating first over parallel hyperplanes, we see that the right side is at most
$$
\int_{\substack{u\in D\\0<\langle u,\theta\rangle-a_\theta<\frac{d_\theta}{2}}}\frac{|f(u)|}{\langle u,\theta\rangle-a_\theta}du.
$$
For fixed $\theta$ put $v(u)=\frac{\langle u,\theta\rangle-a_\theta}{d_\theta}$. When $0<v(u)<\frac{1}{2}$ we have $\int_{v(u)}^1\frac{ds}{s^2}=\frac{1}{v(u)}-1\geq\frac{1}{2v(u)}$, and therefore
$$
\frac{1}{\langle u,\theta\rangle-a_\theta}=\frac{1}{d_\theta v(u)}\leq\frac{2}{d_*}\int_{v(u)}^1\frac{ds}{s^2}.
$$
By Tonelli's theorem,
$$
\int_{\substack{u\in D\\0<v(u)<\frac{1}{2}}}|f(u)|\left(\int_{v(u)}^1\frac{ds}{s^2}\right)du=\int_0^1\frac{1}{s^2}\int_{\substack{u\in D\\0<v(u)<\min\{s,\frac{1}{2}\}}}|f(u)|du ds,
$$
and since the points $u$ with $0<v(u)<s$ belong to $D_s$ by \eqref{eq_layer}, the endpoint contribution is at most
$$
\frac{2}{d_*}\int_0^1\frac{1}{s^2}\int_{D_s}|f(u)|du ds.
$$
The endpoint $b_\theta$ is treated in the same way, and integrating over the sphere only introduces its finite surface measure.

For the remainder, \eqref{rftab} and \eqref{eq_radon_gradient} imply that
$$
\int_{\mathbb S^{d-1}}\int_2^\infty|F(\theta,t)|dt d\theta\leq C\int_{\mathbb S^{d-1}}\left(d_\theta+\frac{1}{d_\theta}\right)\int_{a_\theta}^{b_\theta}\left|\frac{\partial}{\partial s}Rf(\theta,s)\right|ds d\theta\leq C_D\int_D|\triangledown f(u)|du,
$$
where the last constant is finite because of the bounds on $d_\theta$. Combining the endpoint and remainder estimates with \eqref{eq_polar_weight} proves \eqref{estintft}.
\end{proof}

The Sobolev assumption in Theorem~\ref{avft} allows a stronger conclusion by a classical Fourier inequality. The work of Bourgain and of Pe{\l}czy{\'n}ski and Wojciechowski gives, for $d\geq2$ and $F\in W^{1,1}(\mathbb R^d)$,
\begin{equation}\label{eq_classical_sobolev_fourier}
\int_{\mathbb R^d}\frac{|\widehat F(x)|}{|x|^{d-1}}dx
\leq C_d\int_{\mathbb R^d}|\triangledown F(u)|du.
\end{equation}
We refer to \cite{PW} and to the statement in \cite[Theorem 6.1]{Kolyada2019}. Taking $F$ to be the extension of $f$ by zero, whose weak gradient was identified in the proof above, gives
$$
\int_{\mathbb R^d}\frac{|\widehat f(x)|}{|x|^{d-1}}dx
\leq C_d\int_D|\triangledown f(u)|du<\infty.
$$
Thus weighted Fourier integrability under the hypotheses of Theorem~\ref{avft} already follows from this known estimate. It requires neither convexity of the Radon transforms nor finiteness of the boundary layer integral. The proof of \eqref{estintft} gives a direct estimate from the directional representation and retains the boundary layer contributions arising from its endpoint terms.

The weight in \eqref{estintft} cannot simply be omitted, and the example that shows this also shows that the hypotheses allow nonzero subharmonic functions on bounded domains. Let $D$ be the unit disk in $\mathbb R^2$ and take $f(u)=-\sqrt{1-|u|^2}$. This function is continuous on $\overline D$ and vanishes on $\partial D$, and a direct calculation in polar coordinates shows that $\Delta f(u)=\frac{2-|u|^2}{(1-|u|^2)^{\frac{3}{2}}}>0$ and $\int_D|\triangledown f(u)|du=\frac{\pi^2}{2}$. Its extension by zero belongs to $W^{1,1}(\mathbb R^2)$. To see this, integrate by parts against a test function on the disk of radius $r<1$ and let $r$ increase to one. The boundary term tends to zero because $f$ tends to zero uniformly on the circle of radius $r$, and the gradient is integrable, so the weak gradient of the extension is the extension by zero of $\triangledown f$. The chord integral is the negative of the area of a half disk of radius $\sqrt{1-t^2}$, so $Rf(\theta,t)=-\frac{\pi}{2}(1-t^2)$ for $|t|<1$, which is convex. Here $a_\theta=-1$, $b_\theta=1$ and $d_\theta=1$, so $D_t=\{u\in D\mid1-t<|u|<1\}$ and $\int_{D_t}|f(u)|du=\frac{2\pi}{3}(2t-t^2)^{\frac{3}{2}}$. Near $t=0$ this is bounded by a constant times $t^{\frac{3}{2}}$, so its product with $t^{-2}$ is bounded by a constant times $t^{-\frac{1}{2}}$, which is integrable, and away from zero there is no singularity. The boundary layer integral in \eqref{estintft} is therefore finite. To compute the Fourier transform we integrate the polynomial $Rf(\theta,\cdot)$, and since it is even, for $r>0$
$$
\widehat f(r\theta)=-\frac{\pi}{2}\int_{-1}^1(1-t^2)e^{-irt}dt=-\pi\int_0^1(1-t^2)\cos(rt)dt.
$$
Integrating by parts twice, we find that
$$
\int_0^1(1-t^2)\cos(rt)dt=\frac{2}{r}\int_0^1t\sin(rt)dt=-\frac{2\cos r}{r^2}+\frac{2\sin r}{r^3},
$$
and therefore
$$
\widehat f(x)=\frac{2\pi(r\cos r-\sin r)}{r^3}
$$
for $r=|x|>0$. The weighted integral in \eqref{estintft} converges, but
$$
\int_{|x|>2}|\widehat f(x)|dx=4\pi^2\int_2^\infty\frac{|r\cos r-\sin r|}{r^2}dr=\infty.
$$
Indeed, the last integrand is at least $\frac{|\cos r|}{r}-\frac{1}{r^2}$, the first term is at least $\frac{1}{2r}$ on each interval $[2\pi k,2\pi k+\frac{\pi}{3}]$, and its integrals over these intervals dominate a constant multiple of the harmonic series, whereas $\frac{1}{r^2}$ is integrable. Thus even zero boundary values and convexity of every Radon transform do not imply integrability of the Fourier transform under our assumptions.

\subsection{Multidimensional convex functions}

A finite real-valued function $f$ on an open convex set $D\subset\mathbb R^d$ is convex if
\begin{equation}\label{clconv}
f((1-\lambda)x+\lambda y)\leq(1-\lambda)f(x)+\lambda f(y)
\end{equation}
for all $x,y\in D$ and $\lambda\in[0,1]$. Such a function is locally Lipschitz, and hence belongs to $W^{1,1}_{\mathrm{loc}}(D)$, but this local regularity says nothing about integrability at the boundary or at infinity. We refer to \cite[Chapter II]{Hormander} for the basic theory. The stronger theorem on almost everywhere second differentiability, due to Busemann, Feller and Aleksandrov, is discussed in \cite{BCP}, and related multidimensional Fourier estimates are treated in \cite[Chapter 7]{L2019}.

Every convex function is subharmonic. For a twice continuously differentiable function, convexity means that its Hessian is positive semidefinite, so its Laplacian, the trace of the Hessian, is nonnegative. For a general convex function, convolution with a nonnegative smooth approximate identity preserves convexity on each smaller convex domain, and the convolutions converge locally uniformly to $f$. If $\psi\in C_c^\infty(D)$ is nonnegative, integration by parts for these smooth approximations and passage to the limit show that $\int_Df\Delta\psi\geq0$.

\begin{remark}\label{avftc}
Let $D$ be a bounded admissible domain, and let $f$ be convex on $D$. If its extension by zero belongs to $W^{1,1}(\mathbb R^d)$ and $Rf(\theta,\cdot)$ has a convex representative on $(a_\theta,b_\theta)$ for almost every $\theta$, then \eqref{estintft} holds. Indeed, the restriction of the zero extension to $D$ belongs to $W^{1,1}(D)$, so Theorem~\ref{avft} applies.
\end{remark}

The hypothesis on $Rf$ does not follow from convexity of $f$. The function $f(x,y)=x^2+y^2-1$ on the unit disk is convex and vanishes on the boundary, but its Radon transform is not convex, as computed in Section~\ref{sec_examples}. There are also nonzero examples satisfying all the hypotheses of Remark~\ref{avftc}. The function $f(u)=-\sqrt{1-|u|^2}$ on the unit disk, considered above, is convex, since for $u$ in the disk and $v\in\mathbb R^2$
$$
\sum_{j,k=1}^2v_jv_k\partial_j\partial_kf(u)
=\frac{|v|^2}{\sqrt{1-|u|^2}}+\frac{\langle u,v\rangle^2}{(1-|u|^2)^{\frac{3}{2}}}\geq0.
$$
Its Radon transforms are convex and its weighted Fourier integral is finite, while its unweighted Fourier integral diverges.

On a slab, convexity imposes a further restriction. Fix $d\geq2$, $\theta\in\mathbb S^{d-1}$ and finite numbers $a<b$. Suppose that $f$ is convex on
$$
D=\{u\in\mathbb R^d\mid a<\langle u,\theta\rangle<b\},
$$
and that $\int_J\int_{\theta^\perp}|f(t\theta+y)|dy dt<\infty$ for every compact interval $J\subset(a,b)$. Then $f=0$ on $D$.

To prove this, first let $h$ be a finite convex function in $L^1(\mathbb R^m)$, where $m\geq1$. Let $B_R(y)$ be the ball of radius $R$ centered at $y$, and let $|B_R|$ denote its volume. By Jensen's inequality,
$$
h(y)\leq\frac{1}{|B_R|}\int_{B_R(y)}h(z)dz
\leq\frac{\|h\|_{L^1(\mathbb R^m)}}{|B_R|},
$$
and letting $R\to\infty$ shows that $h\leq0$. The restriction of $h$ to any line is a convex function on $\mathbb R$ bounded above, so it is constant. Thus $h$ is constant on $\mathbb R^m$, and integrability forces $h=0$. Fubini's theorem now shows that, for almost every $t\in(a,b)$, the convex function $y\mapsto f(t\theta+y)$ belongs to $L^1(\theta^\perp)$ and therefore vanishes. Since $f$ is continuous on $D$, it vanishes everywhere. Thus the integrability assumptions of Theorem~\ref{corsbh} allow no nonzero convex function on the slab.

\subsection{Concluding remarks}

Recall that a function $f$ is called superharmonic if $-f$ is subharmonic. Of course, the results obtained are
true for superharmonic functions.

By the maximum principle, a subharmonic (seperharmonic) function vanishing on the boundary of an admissible
domain must be non-positive (non-negative).

Observe that in this paper, the one-dimensional and multidimensional directional results on the Fourier transform
are of linear nature. Therefore, they hold true for the so-called quasi-harmonic functions, that is, the functions
that are represented as the difference of two convex functions. Correspondingly, such results can be applied for
functions that can be represented as the difference of two subharmonic (superharmonic) functions.

\end{document}